\documentclass[14pt]{extarticle}
\usepackage[cp1251]{inputenc}
\usepackage{amsmath}
\usepackage{amssymb}

\usepackage[pdftex]{graphicx}

\newcommand{\supp}{\operatorname{supp}}
\newcommand{\Var}{\operatorname{Var}}

\newtheorem{theorem*}{Теорема}
\newtheorem{theorem}{Theorem}
\newtheorem{lem}{Lemma}
\newtheorem{predl}{Proposition}

\numberwithin{equation}{section}

\renewcommand{\hat}{\widehat}

\newcommand{\C}{\mathbb{C}}

\renewcommand{\L}{\Lambda}
\renewcommand{\l}{\lambda}
\renewcommand{\P}{\mathbb{P}}
\newcommand{\E}{\mathbb{E}}

\newcommand{\e}{\varepsilon}

\renewcommand{\a}{\alpha}
\renewcommand{\b}{\beta}
\renewcommand{\d}{\delta}
\newcommand{\g}{\gamma}
\newcommand{\N}{\mathbb{N}}
\newcommand{\Z}{\mathbb{Z}}
\renewcommand{\leq}{\leqslant}
\renewcommand{\geq}{\geqslant}

\begin{document}

	\begin{center}
		\textbf{LOWER BOUNDS FOR THE WIENER NORM IN $\Z_p^d$}
	\end{center}


\begin{center}
	MIKHAIL GABDULLIN
\end{center}
	
	
	\begin{abstract}
		We obtain lower bounds for the $l_1$-norm of the Fourier transform of functions on $\Z_p^d$.
	\end{abstract} 	
	


\bigskip

	\section{Introduction}
	Let $B\subset\Z$ be a finite set of integers, $|B|\geq2$, and $e(x):=\exp(2\pi ix)$. A famous conjecture of Littlewood asserts that 
	$$
	\int_{0}^{1}\left|\sum_{b\in B}e(bx)\right|dx \gg \log|B|.
	$$
	This inequality was first proved by Konyagin \cite{Kon} in 1981. In the slightly later work \cite{MPS} a more general result was obtained: if $B=\{b_1<\ldots<b_n\}$ and $c(b_j)\in\C$ are arbitrary complex numbers then
	\begin{equation}\label{Littlewood} 
	\int_{0}^{1}\left|\sum_{b\in B}c(b)e(bx)\right|dx \gg \sum_{j=1}^n\frac{|c(b_j)|}{j}. 
	\end{equation} 
	It is well-known (see, for instance, \cite{Teml} p.19) that the $L_1$-norm of the Dirichlet kernel $D_n(x)=\sum_{b\in B}e(bx)=\sum_{k=-n}^ne(kx)$ (which corresponds to the case $B=[-n,n]\cap\Z$) is of the order of $\log|B|$. Thus both of the mentioned results are sharp up to constants.

	In the work of Green and Konyagin \cite{GK} and in several subsequent papers
	 \cite{Sand1}, \cite{KS small and big}, \cite{KS medium}, \cite{Sch}, \cite{Sand2} a discrete analog of the Littlewood conjecture (for the case of group $\Z_p$) has been studied. We need some basic definitions (see, for example, \cite{TV}, Chapter 4). Let $G$ be a finite abelian group. A character of the group $G$ is a homomorphism $\g\colon G\to S^1$, where $S^1=\{z\in\C:|z|=1\}$. Denote by $\hat{G}$ the dual group of the group $G$, that is, the group of characters of $G$ with pointwise multiplication as the group operation. 
	 It is well-known that $G$ and $\hat{G}$ 
	 are isomorphic whenever $G$ is abelian and finite; we identify them. For any function $f\colon G\to\C$, its Fourier transform is the function $\hat{f}\colon G\to\C$ defined by
	$$\hat{f}(\g)=|G|^{-1}\sum_{x\in G}f(x)\overline{\g(x)}; 
	$$
the Wiener norm of $f$ is the $l_1$-norm of its Fourier transform, that is,
	$$
	\|\hat{f}\|_1=\sum_{\g\in \hat{G}}|\hat{f}(\g)| 
	$$
	(we write $\|g\|_q:=(\sum_{x\in \supp g}|g(x)|^q )^{1/q}$ for any $q>0$). It is easy to see that
	\begin{equation*}
	\hat{fg}(\xi)=\sum_{\eta\in \hat{G}}\hat{f}(\xi-\eta)\hat{g}(\eta),  
	\end{equation*} 
	which immediately implies 
	\begin{equation}\label{banach} 
	\|\hat{fg}\|_1\leq\|\hat{f}\|_1\|\hat{g}\|_1.  
	\end{equation}
	Thus the algebra of functions on $G$ endowed with the Wiener norm forms a Banach algebra. 
	
	For any set $A\subseteq G$ we denote by $A(\cdot)$ its indicator function, that is, $A(x)=1$ if $x\in A$ and $A(x)=0$ otherwise.  For $u\in\Z_p$, we set  $e_p(u)=e^{2\pi iu/p}$.
	
	In the case $G=\Z_p$ the characters are the functions $\g_{\xi}(x)=e_p(\xi x)$, where $\xi\in\Z_p$,
	and the Wiener norm of the indicator function of a set $A\subseteq\Z_p$ is the discretization of the $L_1$-norm of the polynomial $T_A(y)=\sum_{x\in A}e(xy)$  :
	$$\|\hat{A}\|_1=\frac1p\sum_{\xi \in\Z_p}\left|\sum_{x\in A}\exp\left(\frac{2\pi i \xi x}{p}\right)\right|=\frac1p\sum_{\xi=1}^p\left|T_A\left(\frac{\xi}{p}\right)\right|. 
	$$
	Therefore the problem of finding lower bounds for the Wiener norm of subsets of $\Z_p$ can be considered as a discrete version of the Littlewood conjecture.
	
	It is easy to see that
	$$
	\|\hat{A}\|_1=\|\hat{(\Z_p\setminus A)}\|_1+\frac{2|A|}{p}-1=\|\hat{(\Z_p\setminus A)}\|_1+O(1),
	$$
	and hence it is enough to consider the case $|A|<p/2$. It is highly believed that for all such $A$ the bound 
	\begin{equation}\label{hypoth}
	\|\hat{A}\|_1\gg \log |A|
	\end{equation}
	holds, which is analogous to the corresponding bound (\ref{Littlewood}) in the continious setting.
	
	Let us discuss the known results in this direction. 
	
	Firstly, if $A$ is an arithmetic progression (and $|A|< p/2$) then (\ref{hypoth})
	 holds: we can apply the Marcienkiewicz theorem on discretizations of integral norms for trigonometric polynomials (see, for instance, Theorem 1.10 for the $L_1$-norm in \cite{Teml}) and get
	\begin{equation*}
	\|\hat{A}\|_1\asymp \log |A|.
	\end{equation*}

	Further, the inequality (\ref{hypoth}) was proved by Konyagin and Shkredov for subsets of small size.
	
	\bigskip
	
	\textbf{Theorem A \cite{KS small and big}.} \textit{Let $A\subset \Z_p$ and
		$$|\hat{A}|_1\ll \exp((\log p/\log\log p)^{1/3}).
		$$ 
		Then the bound (\ref{hypoth}) holds.}
	
	Also they obtained the following bounds.
	
	\bigskip
	
	\textbf{Theorem B \cite{KS medium}.} 
	\textit{Let $A\subset\Z_p$ be such that $\exp((\log p/\log\log p)^{1/3})\leq |A|\leq p/3$ and set $\d=|A|/p$. Then}
	$$\|\hat{A}\|_1\gg (\log \d^{-1})^{1/3}(\log\log\d^{-1})^{-1-o(1)}, \quad \d\to0.
	$$		
	
	\bigskip
	
Besides, it was shown in the work \cite{KS small and big} that the results of Sanders \cite{Sand2} imply bounds for dense and near-dense subsets.

	\bigskip
	
	\textbf{Theorem C \cite{KS small and big}.} 
	\textit{Let $A\subset\Z_p$ and $\d=|A|/p<1/2$. Then we have
		$$\|\hat{A}\|_1\gg \d^{3/2}(\log p)^{1/2-o(1)}
		$$
		for $\d\geq (\log p)^{-1/4}(\log\log p)^{1/2}$, and
		$$\|\hat{A}\|_1\gg \d^{1/2}(\log p)^{1/4-o(1)}
		$$
		for} $\d< (\log p)^{-1/4}(\log\log p)^{1/2}$.
	
	\bigskip
	
Theorems B and C give bounds of the type $\|\hat{A}\|_1\gg (\log p)^c$, $c>0$, for dense and sparse subsets, but for example in the case $\d \asymp (\log p)^{-1}$ (where we can only use Theorem B) we have just a lower bound of order $(\log\log p)^{1/3}$. Schoen \cite{Sch} proved that $\|\hat{A}\|_1\gg (\log|A|)^{1/16-o(1)}$ for any $|A|$ (but, as always, under the condition $|A|<p/2$). The following result of Sanders gives a stronger bound $\|\hat{A}\|_1\gg(\log|A|)^{1/4-o(1)}$. 
	
	\bigskip
	
	\textbf{Theorem D \cite{Sand2}.}
	\textit{Let $G$ be a finite abelian group and denote by $\mathcal{W}(G)$ the set of cosets of subgroups of $G$. Then for any function $f\colon G\to\Z$ with $\|\hat{f}\|_1\leq K$ there exists some $z\colon \mathcal{W}(G)\to\Z$ such that
	}
	$$f=\sum_{W\in\mathcal{W}(G)}z(W)1_W
	$$		
	\textit{and $\sum_{W\in\mathcal{W}(G)}|z(W)|\leq\exp(K^{4+o(1)})$.}

	\bigskip

	In this paper we have two purposes. The first one is to generalize Theorems A and B for the functions obeying $|f(x)|\geq1$ for all $x\in\supp f$ instead of indicator functions.
	

	\begin{theorem}\label{th1} 
		Let $f\colon \Z_p\to\C$ be such that $|f(x)|\geq 1$ for all $x\in S:=\supp f$, and set $M:=\max_{x\in\Z_p}|f(x)|$. Then we have
		$$\|\hat{f}\|_1\geq M
		$$
		and
		$$\|\hat{f}\|_1\gg \min\left(\log|S|, \left(\frac{\log p}{(\log\log p)(\log |S|) }\right)^{1/2} \right).
		$$
		In particular, we have $\|\hat{f}\|_1\gg\log|S|$ whenever $|S|\leq \exp((\log p/\log\log p)^{1/3})$. 
	\end{theorem}

	\begin{theorem}\label{th2} 
		Let $f\colon \Z_p\to\C$ be such that $|f(x)|\geq1$ for all $x\in S:=\supp f$, and set $M:=\max_{x\in\Z_p}|f(x)|$. Suppose that $\exp((\log p/\log\log p)^{1/3})\leq|S|\leq p/3$; then 
		$$\|\hat{f}\|_1\geq M
		$$
		and
		$$\|\hat{f}\|_1\gg  (\log (p/|S|))^{1/3}(\log\log(p/|S|))^{-1-o(1)}, \quad p/|S|\to\infty.
		$$
	\end{theorem}
	

The proofs of Theorems 1 and 2 are just technical modifications of the proofs of Theorems A and B.

	The other purpose of this paper is to transfer one-dimensional results to the multidimensional case. We have two results in this direction.

	\begin{theorem}\label{th3} Let $E\subseteq\C$ be arbitrary set and $C>0$ be large absolute constant. Suppose that for any function $h\colon \Z_p\to E\cup\{0\}$ a bound
		$$\|\hat{h}\|_1\geq F(p,\d),
		$$	
		holds, where $\d=|\supp h|p^{-1}$. Then for any function $f\colon\Z_p^d\to E\cup\{0\}$ such that $|\supp f|=\d p^d$, $\d\geq Cp^{-1}$ we have
		$$\|\hat{f}\|_1\geq F(p,\d')
		$$	
		for some $\d'=\d+O(\d^{1/2}p^{-1/2})$, with the implied constant being absolute.
	\end{theorem}
	
	\bigskip
	
	\textit{\textbf{Remark.} Without loss of generality we can assume that $F(p,\d)$ is equal to $+\infty$ if $p\d$ is not an integer. This makes no change in the assumption of the theorem  and excludes trivial cases which are technically possible since we do not know the exact value of $\d'$.}
	
	\bigskip
	
	In particular, one can estimate the Wiener norm of large subsets of $\Z_p^d$ using Theorem \ref{th3} and Theorem C: for any $A\subset\Z_p^d$ with $|A|\asymp p^d$ (and $|A|<p^d/2$) we have $\|\hat{A}\|_1\gg (\log p)^{1/2-o(1)}$. 
	
	It is convenient for us to formulate Theorem \ref{th3} in such a ``conditional'' form, in order to not depend on the best known results in one-dimensional case. Moreover, as was mentioned before, for various classes of functions various bounds are known: for instance, in \cite{Sand2} integer-valued functions are studied, whereas in \cite{KS small and big} and \cite{KS medium} indicators functions are learned, and this paper we consider functions which take values in the set $\{z\in \C: |z|\geq1 \}$. 
	
	In the same conditional form we formulate a result on the Wiener norm of small subsets of $\Z_p^d$. 
	
	\begin{theorem}\label{th4} 
		Let $E\subseteq\C$ be any set invariant under rotation (that is, $e^{i\varphi}E=E$ for all $\varphi\in\mathbb{R}$). Suppose that for any function $h\colon\Z_p\to E\cup\{0\}$ with $|\supp h|=\d p<(2p)^{1/2}$ a bound
		$$\|\hat{h}\|_1\geq F(p,\d)
		$$
		holds. Then for any function $f\colon\Z_p^d\to E\cup\{0\}$ with $|\supp f|=\d p<(2p)^{1/2}$ we have 
		$$\|\hat{f}\|_1\geq F(p,\d).
		$$ 
	\end{theorem} 
	
	In particular, using Theorems 1 and 4 (in the cases $M=1$ and $E=\{z\in\C: |z|=1\}$ respectively) one can estimate the Wiener norm of small subsets of $\Z_p^d$: for $A\subset\Z_p^d$ such that $|A|\leq \exp((\log p/\log\log p)^{1/3})$ we get the tight bound $\|\hat{A}\|_1\gg\log|A|$. 

	For $d\geq2$ it is, generally speaking, we cannot estimate $\|\hat{A}\|_1$ from below by a function which is increasing on $|A|$ since (as easy to see) we have $\|\hat{V}\|_1=1$ for any subroup $V\subseteq \Z_p^d$. Nevertheless, Theorem D gives non-trivial lower bounds for the cases where a set $A$ is far from subgroups: for example, if $G=\Z_p^d$, a number $\eta\in(0,1)$ is fixed and $|A|\asymp p^{k+\eta}$ where $k\in\{0,1,\ldots,d-1\}$ then we get $\|\hat{A}\|_1\gg (\log p)^{1/4-o(1)}$. Theorems 3 and 4 improve this bounds for small and large subsets.
	

	To prove Theorem 3 we use the probabilistic method for finding a line $l\subset \Z_p^d$ such that $p^{-1}|l\cap \supp f|$  is close to $p^{-d}|\supp f|$ and then show that the Wiener norm of a function $f$ in $\Z_p^d$ is no less than the Wiener norm of its restriction to any line (see the precise definition of restriction in Section \ref{3}). To prove Theorem 4, we show that the image of a set $\supp f$ under an appropriate non-degenerate linear map has the property that the first coordinates of its elements are distinct.
	

	\section{Proofs of Theorems 1 and 2}
	
	We follow the works \cite{KS small and big} and \cite{KS medium}.

	\bigskip
	
	Firstly, let us note that the Fourier inversion formula 
	$$f(x)=\sum_{\xi\in\Z_p}\hat{f}(\xi)e_p(\xi x) 
	$$	
	 immediately implies
	\begin{equation}\label{inverse} 
	M=\max_{x\in\Z_p}|f(x)|\leq \sum_{\xi\in\Z_p}|\hat{f}(\xi)|=\|\hat{f}\|_1,
	\end{equation}
	and thus we get the first inequalities in Theorems 1 and 2. 
	Moreover, we see that while proving lower bounds for $\|\hat{f}\|_1$ we can assume that $f$ does not take large values.

	\subsection{Proof of Theorem 1}
	
	
	Let $G$ be a finite abelian group. For any function $g\colon G\to\C$ define the quantity 
	$$T_k(g):=\sum_{x_1+\ldots+x_k=x_1'+\ldots+x_k'}g(x_1)\ldots g(x_k)\overline{g(x_1')}\ldots\overline{g(x_k')}.
	$$
	It is easy to show that
	\begin{equation}\label{T_k}
	T_k(g)=|G|^{2k-1}\sum_{\g}|\hat{g}(\g)|^{2k}. 
	\end{equation}
	Further, a set $\L=\{\l_1,\ldots,\l_{|\L|}\}\subseteq G$ is said to be dissociated if the equality
	$$\sum_{\l\in\L}\e_{\l}\l =0
	$$
	with $\e_{\l}\in\{-1,0,1\}$ implies that $\e_{\l}$ are all equal to zero.
	Let $S\subseteq G$; the size of maximal dissociated subset of the set $S$ is called additive dimension of the set $S$ and denoted by $\dim S$.

	\begin{lem}\label{T_k(g)}
	Let $G$  be a finite abelian group and let $f\colon G\to\C$ be a function with
	 $\|\hat{f}\|_1\leq K$. Let, further, a set $Q\subseteq S$ have the property that $|f(x)|\geq L$ for $x\in Q$, ans set $g(x)=Q(x)f(x)$. Then for all $k\in\N$ we have
		$$T_k(g)\geq \frac{|Q|^{2k}L^{4k}}{\|f\|_2^2K^{2k-2}}\,.
		$$	
	\end{lem} 
	
	\textit{Proof.} Applying the lower bound for $|f(x)|$, Plancherel's theorem, and H\"older's inequality (with exponents $2k$ and $2k/(2k-1)$), we get

	\begin{multline*}
	|Q|^{2k}L^{4k}\leq \left(\sum_xg(x)\overline{f(x)}\right)^{2k}\leq\left(|G|\sum_\g|\hat{g}(\g)||\hat{f}(\g)|\right)^{2k} \leq \\
	\leq |G|^{2k}\sum_\g |\hat{g}(\g)|^{2k}\left(\sum_\g|\hat{f}(\g)|^{2k/(2k-1)}\right)^{2k-1}.
	\end{multline*}
	Writing $|\hat{f}(\g)|^{2k/(2k-1)}$ as $|\hat{f}(\g)|^{(2k-2)/(2k-1)}|\hat{f}(\g)|^{2/(2k-1)}$ and using H\"older inequality (with exponents $(2k-1)/(2k-2)$ and $2k-1$), (\ref{T_k}) and upper bound  for $\|\hat{f}\|_1$, we have
	\begin{multline*}
	|Q|^{2k}L^{4k} \leq |G|T_k(g)\left(\sum_\g|\hat{f}(\g)|\right)^{2k-2}\sum_\g|\hat{f}(\g)|^2 \leq \\
	T_k(g)K^{2k-2}\sum_{x}|f(x)|^2 = T_k(g)K^{2k-2}\|f\|_2^2.
	\end{multline*}
	This concludes the proof. \quad $\Box$

	\bigskip

	\begin{predl}
		Let $G$ be a finite abelian group and a function $f\colon G\to\C$ obey  $|f(x)|\geq1$ for $x\in S=\supp f$ and $\|\hat{f}\|_1\leq K\leq \|f\|_2$. Then
		$$\dim S \ll K^2\left(1+\log \frac{\|f\|_2}{K}\right).
		$$
	\end{predl}
	
	\textit{Proof.} Let $\Lambda$ be a maximal dissociated subset of the set $S$, and
	$$\Lambda_j=\{x\in \Lambda : 2^{j-1}\leq |f(x)|<2^j\}, \qquad j=1,2,\ldots.
	$$
	We apply the previous lemma with $Q=\Lambda_j$ and obtain (here $g_j(x)=\Lambda_j(x)f(x)$)
	$$\frac{|\L_j|^{2k}2^{4(j-1)k}}{\|f\|_2^2K^{2k-2}} \leq T_k(g_j)\leq T_k(|g_j|)\leq 2^{2jk}T_k(\Lambda_j).
	$$
	Well-known Rudin's inequality claims that for any dissociated set $\Lambda$ the bound $T_k(\Lambda)\leq (Ck)^k|\Lambda|^k$ holds with an absolute constant $C>0$ (see, for instance, \cite{TV}: it follows from the bound (4.34) and Definition 4.26 $\L(p)$-constant). The sets $\Lambda_j$ are dissociated as subsets of a dissociated $\Lambda$; thus
	$$|\L_j|^k\leq (2^{-2j+4}Ck)^k\|f\|_2^2K^{2k-2},
	$$
	and hence
	$$|\Lambda_j|\ll 2^{-2j}K^2\exp\left(\log k+\frac2k\log\frac{\|f\|_2}{K} \right).
	$$
	Setting $k:=[2\log \frac{\|f\|_2}{K}]+1$, we see that
	$$|\Lambda_j|\ll 2^{-2j}K^2\left(1+\log \frac{\|f\|_2}{K}\right).
	$$
	We are done by taking the sum over all $j\geq1$ and recalling that $|\Lambda|=\dim S$. $\Box$

	\bigskip 
	
	\textit{\textbf{Remark.} By Cauchy-Schwarz and the Parseval identity we see that 
		$$\|\hat{f}\|_1\leq p^{1/2}\left(\sum_{\g}|\hat{f}(\g)|^2\right)^{1/2}=\left(\sum_{x}|f(x)|^2\right)^{1/2}= \|f\|_2,
		$$
	and hence the restriction $K\leq \|f\|_2$ in Proposition 1 is not essential.}
	
	\bigskip
	
	Now we are ready to prove Theorem 1. The inequality (\ref{inverse}) gives the first inequality in the theorem; moreover, it follows that while proving the second one we can restrict our attention to the case $M=\max_{x\in\Z_p}|f(x)|\leq \log|S|$.

	Set $K:=\|\hat{f}\|_1$. Let $c>0$ be a small absolute constant. We can assume that $K^2\log|S| \leq c\frac{\log p}{\log \log p}$ since otherwise we are done. Let $\Lambda$ be a maximal dissociated subset of $S$. Since $M\leq \log|S|$, we have $\|f\|_2\leq |S|^{1/2}M\leq |S|^{1/2}\log|S|$, and, hence, $\log\frac{\|f\|_2}{K}\ll \log|S|$. Then by Proposition 1 we get
	\begin{equation}\label{dim} 
	d:=\dim S=|\L| \leq \frac{\log p}{\log\log p} \, .
	\end{equation}
	For any $x\in\Z_p$ we set $|x|=\min\{|z| : z\in\Z,\, z\equiv x\pmod{p} \}$ . By Dirichlet's approximation theorem there exists a positive integer $q<p$ such that for all  $\lambda\in\Lambda$
	\begin{equation}\label{Dirichlet} 
	|q\lambda|\leq p^{1-1/d}\leq \frac{p}{\log p}. 
	\end{equation}
	
	Now take an arbitrary element $x\in S$; by the definition of the set $\L$ we have $x=\sum_{\l\in\L}\e_{\l}\l$ for some numbers $\e_{\l}\in\{-1,0,1\}$; by   (\ref{dim}) and (\ref{Dirichlet}) we get 
	$$|qx|=|q\sum_{\l\in\L}\e_{\l}\l|\leq\sum_{\l\in\L}|q\l|\leq \frac{p}{\log\log p}\leq \frac{p}{3} 
	$$
	for sufficiently large $p$. Consider the function $F(x):=f(qx)$; its support is the set $B:=qS\subseteq [-p/3,p/3]$ and
	$$\|\hat{f}\|_1=\|\hat{F}\|_1=\frac1p\sum_{\xi=1}^p\left|\sum_{x\in B}F(x)e_p(\xi x)\right|
	$$

	Let $B=\{b_1<\ldots b_n\}$. Consider the trigonometrical polynomial $T(y)=\sum_{j=1}^nF(b_j)e(b_jy)$. Using the Marcienkiewicz theorem and the inequalities (\ref{Littlewood}) and $|F(b_j)|\geq1$, we get
	$$\|\hat{f}\|_1=\|\hat{F}\|_1 \gg \int_0^1\left|F(y)\right|dy\gg \sum_{j=1}^n\frac{|F(b_j)|}{j}\gg \log |B|=\log|S|.  
	$$
	This completes the proof.

	\subsection{Proof of Theorem 2}


	It remains to prove the second bound. We need some auxiliary results.
	
	\begin{lem}\label{gg}
		Let $n\in\N$, $B\subset[-2n,2n]$, $|B|\geq2$, $0<\eta<1/2$, $|B\cap[-n,n]|\geq(1-\eta)|B|$. Let, further, $c(b)$ ($b\in B$) be complex numbers with $|c(b)|\geq1$ for all $b\in B\cap[-n,n]$. Then we have
		$$\int_0^1\left|\sum_{b\in B}c(b)e(bu)\right|du \gg \min\left(\log\frac{1}{\eta},\log|B|\right).
		$$ 
	\end{lem}
	
	 This statement appeared in the paper \cite{KS medium} for the case $c(b)=1$ for all $b\in B\cap[-n,n]$, but its proof still works for the more general setting of Lemma \ref{gg}.
	
	\bigskip 
	
	
	
	Repeating arguments from \cite{KS medium} and using Lemma \ref{gg} and the de la Vall\'ee Poussin kernels one can prove the following. 
	
	\begin{lem}\label{gg min}
	 Let $B\subset\Z_p$, $n\leq p/6$ and $\eta\in(0,1/2)$. Suppose that $|B\cap[-2n,2n]|\geq2$ and
		$$|B\cap[-n,n]|\geq(1-\eta)|B\cap[-2n,2n]|.
		$$
		Then for any function $F\colon \Z_p\to\C$ such that $|F(b)|\geq1$ for all 
		$b\in B=\supp\, F$ we have
		$$\|\hat{F}\|_1\gg \min\left(\log\frac{1}{\eta},\log|B\cap[-2n,2n]|\right).
		$$
	\end{lem}

	Now we are going to obtain an analog of the structural lemma from \cite{KS medium}. Set
	$$S_j=\{x\in S : 2^{j-1}\leq|f(x)|<2^j\}, \quad j=1,2,\ldots, [\log_2M]+1.
	$$
	Lemma \ref{T_k(g)} for $k=2$, $K=\|\hat{f}\|_1$, $Q=S_j$ and $g_j(x)=S_j(x)f(x)$ gives us
	$$
	\frac{|S_j|^42^{8(j-1)}}{\|f\|_2^2\|\hat{f}\|_1^2}\leq T_2(g)\leq 2^{4j}T_2(S_j), 
	$$
	and therefore
	$$T_2(S_j)\gg \frac{|S_j|^42^{4j}}{\|f\|_2^2\|\hat{f}\|_1^2}. 
	$$
	Since $T_2(S)=T_2(\sqcup_jS_j)\geq \sum_jT_2(S_j)$, we have
	\begin{equation}\label{T} 
	T_2(S)\gg \frac{\sum_j|S_j|^42^{4j}}{\|f\|_2^2\|\hat{f}\|_1^2}. 
	\end{equation}
	By H\"older's inequality we get
	$$|S|=\sum_{j}|S_j|2^j2^{-j}\leq \left(\sum_j|S_j|^42^{4j}\right)^{1/4}\left(\sum_j2^{-4j/3}\right)^{3/4}\ll \left(\sum_j|S_j|^42^{4j}\right)^{1/4}
	$$
	and (since $j\leq [\log_2M]+1$)
	\begin{multline*} \|f\|_2^2=\sum_{x\in S}|f(x)|^2 \asymp \sum_{j}|S_j|2^{j}2^j \ll \\ \left(\sum_j|S_j|^42^{4j}\right)^{1/4}\left(\sum_{j}2^{4j/3}\right)^{3/4} \ll M\left(\sum_j|S_j|^42^{4j}\right)^{1/4}.
	\end{multline*} 
	The latter two inequalities imply
	$$\sum_j|S_j|^42^{4j} \gg |S|^3\|f\|_2^2M^{-1}.
	$$
	Combining this with (\ref{T}), we find
	\begin{equation}\label{T2}
	T_2(S)\gg \frac{|S|^3}{M\|\hat{f}\|_1^2} 
	\end{equation}
	
	Repeating arguments from \cite{KS medium} concerned with Balog-Szemer\'edi-Gowers theorem and Freiman's theorem (in fact, just replacing $\|\hat{\chi_A}\|_1$ by $O(M^{1/2}\|\hat{f}\|_1)$), we get the following statement.
	
	\begin{lem}\label{structure} 
		Let $\e>0$, $R>R(\e)$, and $d_{\e}=\log^{3+\e}R$; let, further, a function $f\colon \Z_p\to\C$ obey $|f(x)|\geq1$ for all $x\in S=\supp f$, and $M=\max_{x\in\Z_p}|f(x)|$,  $m=\left[d_{\e}p\left(\frac{|S|}{p}\right)^{1/d_{\e}}\right]$. Suppose that $M^{1/2}\|\hat{f}\|_1\leq R$. Then there exist $x_0\in\Z_p$ and $q\in\Z_p^*$ such
		that the set
		$$B=q(S-x_0)	
		$$	
	 has the property
		$$|B\cap[-m,m]|\geq |S|e^{-d_{\e}}.
		$$
	\end{lem}
	
	Finally, we need the lemma from \cite{KS medium} on upper estimates for $T_k(Q)$ for scattered sets $Q$. 
	
	\begin{lem}\label{scat}
		Let $I,k,m,N$ be positive integers and $Q=\sqcup_{i=1}^IQ_i\subseteq \Z$, where  $Q_i\subseteq[-4^im,-\frac{4^i}{2}m)\cup(\frac{4^i}{2}m,4^im]$, $|Q_i|=N$, and $i$ runs over a subset of $\Z$ of cardinality $I$. Then
		$$T_k(Q)\leq 2^{8k}k^kI^kN^{2k-1}.
		$$
	\end{lem}

	\bigskip\bigskip  
	
	Now we are ready to prove Theorem 2. We can assume that $M\leq \log^{1/3}(p/|S|)$. Fix $\e>0$ and suppose that 
	\begin{equation}\label{badassump} 
	M^{1/2}\|\hat{f}\|_1\leq R, \quad R(\e)\leq R \leq M^{1/2}(\log(p/|S|))^{1/3}(\log\log(p/|S|))^{-1-\e}.  
	\end{equation}
	Our aim is to prove that the inequalities (\ref{badassump}) cannot hold provided that $p/|S|$ is large enough depending on $\e$. Since $\e$ is arbitrary, Theorem 2 will follow. 
	
	Firstly, since $\|\hat{f}\|_1\geq M$, (\ref{badassump}) implies  
	\begin{equation}\label{R^{2/3}}
	M\leq R^{2/3};
	\end{equation}
	and in particular, $RM^{-1/2}\geq R^{2/3}$.
	
	Take $x_0,q,m$ and $B$ according Lemma \ref{structure}. Set $F(x)=f(q(x-x_0))$. Then $\|\hat{F}\|_1=\|\hat{f}\|_1$ and hence
	\begin{equation}
	M^{1/2}\|\hat{F}\|_1\leq R.
	\end{equation} 
	Let $l_0$ be the maximal integer such that $2^lm\leq p/3$, and
	$$D_l=\{b\in B: |b|\leq 2^lm \}, \quad 0\leq l\leq l_0, 
	$$
	$$\eta=\exp(-CRM^{-{1/2}}),
	$$
 where absolute constant $C>0$ is large enough (note $\eta\in(0,1/2)$), and
	$$N=[\eta|S|e^{-d_{\e}}].
	$$
	Note that by the lower bound on $|S|$ and the upper bounds on $M$ and $R$ we have $N\geq1$ for sufficiently large $p$. We consider two cases.
	
	\textbf{Case I.} Suppose there exists $l$ such that $|D_l\setminus D_{l-1}|<N$. By Lemma \ref{structure} $|D_l|\geq|D_0|\geq |S|e^{-d_{\e}}$, and therefore $|D_{l-1}|\geq |D_l|-N\geq |D_l|(1-N|D_l|^{-1})\geq(1-\eta)|D_l|$. Using Lemma \ref{gg min} for the set $B$ and $n=2^{l-1}m$, we get 
	$$\|\hat{F}\|_1\gg \min\left(\log\frac{1}{\eta},\log|D_0|\right).
	$$
	Using (\ref{badassump}) and the inequality $|S|\geq\exp((\log/\log\log p)^{1/3})$, we have (if $p/|S|$ is large enough) 
	$$\log|D_0|\geq \log|S|-d_{\e}\gg \log|S|\geq (\log p/\log\log p)^{1/3}\geq \log \frac{1}{\eta}=CRM^{-{1/2}},  
	$$
and hence
	$$\|\hat{f}\|_1=\|\hat{F}\|_1\gg CRM^{-1/2},
	$$
	and we get a contradiction with (\ref{badassump}) provided that $C$ is sufficiently large.
	
	\bigskip
	
	\textbf{Case II.} Suppose $|D_l\setminus D_{l-1}|\geq N$ for all $l=1,\ldots,l_0$. Then for each $l\equiv 0\pmod{2}$ we can choose a subset $S_l\subseteq D_l\setminus D_{l-1}$ with $|S_l|=N$. Set 
	$$Q=\bigsqcup_lS_l;
	$$
	we can apply Lemma \ref{scat} 
	(with $I=[l_0/2]$) to the set $Q$ and sets $S_l$. 
	 We have
	\begin{equation}\label{Tklow} 
	T_k(Q)\leq (2^8k)^kI^kN^{2k-1}.
	\end{equation}
	Now we obtain a lower bound on $T_k(Q)$. Let
	$$Q_j=\{x\in Q: 2^{j-1}<|F(x)|\leq 2^j \}. 
	$$
	Applying Lemma \ref{T_k(g)} to the set $Q_j\subseteq B$, number $K=RM^{-1/2}$, and functions $F(x)$ and $g_j(x)=Q_j(x)F(x)$, we get
	$$\frac{|Q_j|^{2k}2^{4(j-1)k}}{\|F\|_2^2K^{2k-2}}\leq 
	T_k(g_j)\leq T_k(|g_j|)\leq 2^{2jk}T_k(Q_j).  
	$$
	Therefore
	$$T_k(Q_j) \geq \frac{|Q_j|^{2k}2^{2jk-4k}}{\|F\|_2^2K^{2k-2}};
	$$
	summing this over all $j$, we find
	$$T_k(Q)\geq \sum_jT_k(Q_j)\geq 
	\frac{1}{2^{4k}\|F\|_2^2K^{2k-2}}\sum_{j}|Q_j|^{2k}2^{2jk}.  
	$$
	Further, since
	\begin{multline*}
	|Q|^{2k}=\left(\sum_j|Q_j|2^{j}2^{-j}\right)^{2k}\leq \left(\sum_j|Q_j|^{2k}2^{2jk}\right)\left(\sum_j2^{-2jk/(2k-1)}\right)^{(2k-1)}  \\
	\leq \sum_j|Q_j|^{2k}2^{2jk},
	\end{multline*}
	and $|Q|=IN$, we see that
	\begin{equation}\label{Tkup} 
	T_k(Q)\geq \frac{|Q|^{2k}}{2^{4k}\|F\|_2^2K^{2k-2}}=\frac{(IN)^{2k}}{2^{4k}\|F\|_2^2K^{2k-2}} .  
	\end{equation}
	Let us compare the estimates (\ref{Tklow}) and (\ref{Tkup}). It follows that
	$$I\ll K^2k\left(\frac{\|F\|_2^2}{NK^2}\right)^{1/k}\leq K^2k\left(\frac{|S|M^2}{NK^2}\right)^{1/k}.
	$$
	Recalling $K=RM^{-1/2}$ and using the inequality (\ref{R^{2/3}}) we obtain
	\begin{multline} I\ll R^2M^{-1}k\left(\frac{|S|M^3}{NR^2}\right)^{1/k}\leq R^2M^{-1}k\left(\frac{|S|}{N}\right)^{1/k}=\\R^2M^{-1}\exp\left(\log k+\frac1k\log\frac{|S|}{N}\right).
	\end{multline} 
	Set $k=[\log\frac{|S|}{N}]\geq1$; then
	$$I\ll R^2M^{-1}\log\frac{|S|}{N}.
	$$
	Since $|S|N^{-1}\ll \eta^{-1}e^{d_{\e}}\leq \exp(CRM^{-1/2}+\log^{3+\e}R)$ and $RM^{-1/2}\geq R^{2/3}$, we find that $\log\frac{|S|}{N}\ll RM^{-1/2}$; thus
	\begin{equation}\label{I}
	I\ll R^3M^{-3/2}.
	\end{equation}
	But
	$$I\geq l_0/2-1\gg \log(p/m)\gg d_{\e}^{-1}\log(p/|S|)-\log d_{\e}.
	$$
	Recalling $d_{\e}=\log^{3+\e}R$ and taking into account the estimate for $R$ 
	(assumption (\ref{badassump})), we get
	$$I\gg d_{\e}^{-1}\log(p/|S|)\gg \log(p/|S|)(\log\log(p/|S|))^{-3-\e}.
	$$
	Hence, using (\ref{I}) we see that 
	$$R^3 \gg M^{3/2}I\gg M^{3/2}\log(p/|S|)(\log\log(p/|S|))^{-3-\e},
	$$
	which contradicts (\ref{badassump}) whenever $p/|S|$ is large enough. This concludes the proof of Theorem 2.

	\section{Proof of Theorem \ref{th3}}\label{3}
	
	Let $d\geq1$. For $x=(x_1,\ldots,x_d), \xi=(\xi_1,\ldots,\xi_d)\in\Z_p^d$ and $c\in\Z_p$ we set
	$$x\xi=(x,\xi):=x_1\xi_1+\ldots+x_d\xi_d\in\Z_p$$
	and	
	$$
	c\xi=(c\xi_1,\ldots,c\xi_d)\in\Z_p^d.
	$$
	Since the characters of $\Z_p^d$ are the maps
	$$\gamma_{\xi}(x)=e_p(\xi x),
	$$
	where $e_p(u)=\exp(2\pi iu/p)$ and $\xi\in\Z_p^d$, we see that the Fourier transform of a function $f\colon \Z_p^d \to \C$ can be rewritten as
	$$\hat{f}(\xi)=p^{-d}\sum_{x\in\Z_p^d}f(x)e_p(-\xi x),
	$$
	and its Wiener norm is just
	$$\|\hat{f}\|_1=\sum_{\xi\in\Z_p^d}|\hat{f}(\xi)|.
	$$
\bigskip
	In what follows we denote by $\theta$, $\theta'$, $\theta''$ real numbers 
	(which can vary from line to line) bounded in magnitude by $1$.

		Define a line $l\subset\Z_p^d$ to be any set of the form
	\begin{equation}\label{l} 
	l=\{y\in \Z_p^d: y=ub+c \,|\,u\in\Z_p \},
	\end{equation}
	where $c=(c_1,\ldots,c_d)$ and $b=(b_1,\ldots,b_d)$, $b\neq0$, are elements of $\Z_p^d$. A hyperplane is a set $L\subseteq\Z_p^d$ of the form  
	$$L=L_{\eta,u}=\{x\in\Z_p^d: x\eta=u\},
	$$
	where $\eta\in\Z_p^d$, $\eta\neq0$, and $u\in\Z_p$.
	
	We need the following lemmas.
	
	\begin{lem}\label{indstep} 
		Let $A\subset \Z_p^d$ and $|A|=\d p^d$. Then there exists a hyperplane\\ $L\subset \Z_p^d$ such that
		$$|A\cap L|=p^{d-1}(\d+\theta \delta^{1/2} p^{-(d-1)/2}),
		$$
		where $|\theta|\leq1$.
	\end{lem}
	
	\textit{Proof.} We use the probabilistic method. Let us consider all possible ``directions'' $\eta_1,\ldots,\eta_r\in \Z_p^d$, that is, a maximal collection of vectors from $\Z_p^d$ any pair of which is linearly independent. Then each non-zero vector $\xi\in\Z_p^d$ is parallel to exactly one $\eta_j$ (that is, $\xi=\a\eta_j$ for some $\a\in\Z_p$). Thus each direction corresponds to exactly $p-1$ non-zero vectors, and hence $r=\frac{p^d-1}{p-1}=1+p+\ldots+p^{d-1}$.
	
 Now we choose an ordered pair $(\eta,u)\in\L\times\Z_p$ (that is, a hyperplane) uniformly at random and consider the random variable
	$$\xi=|A\cap L_{\eta,u}|=\sum_{x\in A}1(x\eta=u).
	$$
	Firstly,
	$$\E\xi=\frac{1}{rp}\sum_{x\in A}\sum_{\eta\in\L}\sum_{u\in \Z_p}1(x\eta=u)=\frac{1}{rp}\sum_{x\in A}\sum_{\eta\in\L}1=\frac{|A|}{p}.  
	$$
	Further,
	\begin{multline*}
	\E\xi^2=\frac{1}{rp}\sum_{x_1,x_2\in A}\sum_{\eta\in\L}\sum_{u\in \Z_p}1(x_1\eta=x_2\eta=u)=\\\E\xi+\frac{1}{rp}\sum_{x_1,x_2\in A, x_1\neq x_2}\sum_{\eta\in\L}\sum_{u\in \Z_p}1(x_1\eta=x_2\eta=u).
	\end{multline*}
	Now fix $x_1,x_2\in A$, $x_1\neq x_2.$ We say that vectors $a,b\in\Z_p^d$ are orthogonal if $(a,b)=0$. Clearly, $\sum_{u\in \Z_p}1(x_1\eta=x_2\eta=u)=1$ if $(x_1-x_2,\eta)=0$ and $\sum_{u\in \Z_p}1(x_1\eta=x_2\eta=u)=0$ otherwise. Further, there are $p^{d-1}-1$ non-zero vectors $y\in\Z_p^d$ with $(x_1-x_2,y)=0$. Note that all vectors parallel to some $\eta_j$ are simultaneously orthogonal or not orthogonal to the difference $x_1-x_2$. Therefore, the number of $\eta_j\in\L$ which are orthogonal to $x_1-x_2$ is equal to $\frac{p^{d-1}-1}{p-1}$. Recalling that $r=\frac{p^d-1}{p-1}$ we get
	\begin{multline*}
	\E\xi^2=\E\xi+\frac{1}{rp}\sum_{x_1,x_2\in A, x_1\neq x_2}\sum_{\eta\in\L,\eta\bot x_1-x_2}1=\E\xi+\frac{p^{d-1}-1}{rp(p-1)}|A|(|A|-1)=\\
	\E\xi+\frac{1}{p^2}\frac{p^{d-1}-1}{p^{d-1}-1/p}|A|(|A|-1)<\E\xi+\frac{1}{p^2}|A|^2=\E\xi+(\E\xi)^2.
	\end{multline*}
	Then
	$$\Var\xi =\E\xi^2-(\E\xi)^2<\E\xi=\d p^{d-1}. 
	$$
	and
	$$\sigma:=(\Var\xi)^{1/2}=\frac{\d^{1/2}p^{(d-1)/2}}{\l}
	$$
	for some $\l>1$. Thus by Chebyshev's inequality 
	\begin{equation*}
	\P(|\xi-\E\xi|\geq \d^{1/2}p^{(d-1)/2})=\P(|\xi-\E\xi|\geq \l\sigma)\leq \frac{1}{\l^2}<1,
	\end{equation*} 
	and hence there exists a hyperplane $L_{\eta,u}$ such that
	$$|A\cap L_{\eta,u}|=\d p^{d-1}+\theta\d^{1/2}p^{(d-1)/2} 
	$$
	for some $\theta<1$. The claims follows. 
	  \qquad $\Box$
	
	\bigskip

	Now we iterate Lemma \ref{indstep} to find a line $l$ such that $p^{-1}|A\cap l|$ is close to $\d$.  
	
	\begin{lem}\label{to d=1}
		Let $A\subset \Z_p^d$, $|A|=\d p^d$ and $\d\gg p^{-1}$. Then there exists a line $l\subset \Z_p^d$ such that
		$$|A\cap l|=p\left(\d+\theta\d^{1/2}p^{-1/2}+O(\d^{1/2}p^{-1}+p^{-3/2})\right),
		$$
	  for some $|\theta|\leq1$.	
	\end{lem}

	\textit{Proof.} Let $\{e_1,\ldots,e_d\}$ be the standard basis $\Z_p^d$ over $\Z_p$ (in what follows we use only this basis). For $k=1,\ldots,d-1$ we identify the subspaces $\{x\in\Z_p^d: x_{d-k+1}=\ldots=x_{d-1}=x_d=0\}$ with $\Z_p^{d-k}$.
	
	Note if $T\colon\Z_p^d\to\Z_p^d$ is a non-degenerate affine map and $l$ is a line in $\Z_p^d$, then the set $T^{-1}l$ is again a line and $|TA\cap l|=|A\cap T^{-1}l|$. Thus it is enough to prove the lemma for an arbitrary image $A$ with respect to a non-degenerate affine map.
	
	The claim for $d=2$ is just Lemma \ref{indstep}. Let $d\geq3$. First, we show that in this case we can assume  
	\begin{equation}\label{main case Z_p^3}
	\d_{d-3}:=p^{-3}|A\cap \Z_p^3|=\d+\theta\sum_{k=3}^{d-1}p^{-k/2}.
	\end{equation}
	Indeed, for the case $d=3$ we have $p^{-3}|A\cap\Z_p^3|=p^{-3}|A|=\d$. For $d\geq4$ we apply Lemma \ref{indstep} to the set $A$ and find a hyperplane $L_1\subset \Z_p^d$ such that 
	$$p^{-(d-1)}|A\cap L_1|=\d_1,
	$$
	where $\d_1=\d+\theta \d^{1/2}p^{-(d-1)/2}=\d+\theta' p^{-(d-1)/2}$. By the above remark we can think that $L_1=\Z_p^{d-1}$. Using Lemma \ref{indstep} to the set $A\cap L_1\subset \Z_p^{d-1}$, we find a hyperplane $L_2\subset \Z_p^{d-1}$ with 
	$$|A\cap L_2|=\d_1+\theta\d_1^{1/2}p^{-(d-2)/2}=\d_1+\theta'p^{-(d-2)/2}=\d+\theta''\sum_{k=d-2}^{d-1}p^{-k/2}. 
	$$
	Again, without loss of generality we can write $L_2=\Z_p^{d-2}$ and use the lemma to the set $A\cap L_2\subset \Z_p^{d-2}$. Continuing this process, we get (\ref{main case Z_p^3}). 
	
	Thus we can assume that $|A\cap\Z_p^3|=\d_{d-3}p^3$, where $\d_{d-3}=\d+\theta\sum_{k=3}^{d-1}p^{-k/2}=\d+O(p^{-3/2})$. Then we have
	$\d_{d-3}^{1/2}=\d^{1/2}+O(\d^{-1/2}p^{-3/2})$. Applying Lemma \ref{indstep} to the set $A\cap\Z_p^3$ and recalling that $\d\gg p^{-1}$, we find a hyperplane $L\subset\Z_p^3$ such that
	\begin{multline*}
	\d_{d-2}:=p^{-2}|A\cap L|=\d_{d-3}+\theta\d_{d-3}^{1/2}p^{-1}=\\
	\d+O(p^{-3/2}+\d^{1/2}p^{-1}+\d^{-1/2}p^{-5/2})=\d+O(\d^{1/2}p^{-1});
	\end{multline*}
	also we have
	$$\d_{d-2}^{1/2}=\d^{1/2}+O(p^{-1}).
	$$
	As before, we can assume that $L=\Z_p^2$. Applying Lemma \ref{indstep} to the set $A\cap L\subset \Z_p^2$, we find a line $l\subset\Z_p^2$ with 
	\begin{multline*}
	p^{-1}|A\cap l|=\d_{d-2}+\theta\d_{d-2}^{1/2}p^{-1/2}=\d+\theta\d^{1/2}p^{-1/2}+O(\d^{1/2}p^{-1}+p^{-3/2}), 
	\end{multline*}
	as desired. \qquad $\Box$

	

	\bigskip\bigskip
	
Let $l$ be a line defined by (\ref{l}). For a function $f\colon\Z_p^d\to\C$ we define its restriction  $f_l\colon\Z_p\to\C$ to the line  $l$ by
	$$f_l(u)=f(b_1u+c_1,\ldots,b_du+c_d).
	$$
	It is easy to see that the absolute values of the Fourier coefficients of the function $f_l$ (and hence its Wiener norm) do not depend on a parametrization of the line $l$. 
	
	Now we show that the Wiener norm of a function $f$ in $\Z_p^d$ is no less than the Wiener norm of its restriction to any line .

	\begin{lem}\label{f_l}
		Let $f\colon \Z_p^d\to\C$ and $l\subset \Z_p^d$ be a line. Then
		$$\|\hat{f}\|_1 \geq \|\hat{f_l}\|_1.
		$$	
	\end{lem}

\textit{Proof.}
	Let a line $l$ be defined by (\ref{l}) and $l(x)$ be its indicator function. Since $l$ is a shift of one-dimensional subspace, we have $\|\hat{l}\|_1=1$; then the inequality (\ref{banach}) gives us
	\begin{multline*}
	\|\hat{f}\|_1=\|\hat{f}\|_1\|\hat{l}\|_1\geq \|\hat{fl}\|_1=\sum_{\xi\in\Z_p^d}\left|p^{-d}\sum_{x\in l}f(x)e_p(-\xi x)\right|=\\\sum_{\xi\in\Z_p^d}\left|p^{-d}\sum_{u\in \Z_p}f(ub+c)e_p(-\xi(ub+c))\right|=
	p^{-(d-1)}\sum_{\xi\in\Z_p^d}\left|p^{-1}\sum_{u\in \Z_p}f_l(u)e_p(-u\xi b)\right|.
	\end{multline*}
	Since $b\neq0$, for any $\eta\in\Z_p$ there exist exactly $p^{d-1}$ vectors $\xi\in\Z_p^d$ with $\xi b=\eta$. Therefore the latter sum is equal to $$\sum_{\eta\in\Z_p}\left|p^{-1}\sum_{u\in\Z_p}f_l(u)e_p(-u\eta)\right|=\|\hat{f_l}\|_1.
	$$
	The claim follows. $\Box$

Theorem \ref{th3} now follows from Lemmas \ref{to d=1} and \ref{f_l} and the one-dimensional assumption.	
	
	\section{Proof of Theorem \ref{th4}}
	
	 Denote $A=\supp f.$ Firstly, we show that there exists non-degenerate linear map $T\colon\Z_p^d\to\Z_p^d$, $T=\{t_{ij}\}_{i,j=1}^d$, such that the first coordinates of elements of the set $TA$ are distinct. 
	 Let $A=\{a_1,\ldots,a_n\}$ and $a_j=\sum_{i=1}^d\a_{ij}e_i$; then
	$$Ta_j=\sum_{i=1}^d\b_{ij}e_i,
	$$
	where
	$$\b_{ij}=\sum_{k=1}^dt_{ik}\a_{kj}. 
	$$
	We want the sums $\b_{1,j}=\sum_{k=1}^dt_{1,k}\a_{k,j}$, $j=1,\ldots,n$, to be distinct; it is equivalent to say that for all $i<j$ 
	\begin{equation}\label{neq0} 
	\sum_{k=1}^dt_{1,k}(\a_{k,i}-\a_{k,j})\neq0. 
	\end{equation}
	Thus we need to find a vector $\mathbf{t}=(t_{1,1},\ldots,t_{1,d})$ (which will be the first row of the matrix $T$) obeying (\ref{neq0}). There are $|A|(|A|-1)/2$ pairs of elements $(a_i,a_j)\in A^2$, $i<j$; since $|A|<(2p)^{1/2}$, there are less than $p$ such pairs. For each such pair there exist exactly $p^{d-1}$ vectors $\mathbf{t}$ for which (\ref{neq0}) does not hold; since there are $p^d$ vectors in $\Z_p^d$, we can find a vector $\mathbf{t}$ obeying (\ref{neq0}).    
	
	So, there exists non-degenerate linear map $T$ such that the first coordinates of elements of the set $TA$ are distinct. Let 
	$$A'=\{x\in \Z_p: x=\b_{1,j} \mbox{ for some $j$}\}
	$$
	be the ``projection'' of the set $TA$ to the line $x_2=\ldots=x_d=0$. Then each element $a\in TA$ has a representation of the form
	\begin{equation}\label{TA} 
	a=a_x=xe_1+y_2(x)e_2+\ldots+y_d(x)e_d, \quad x\in A'.
	\end{equation} 
	Now consider the function
	$$h(x):=f(T^{-1}x).
	$$
	Then we have
	$$\hat{h}(\xi)=p^{-d}\sum_xf(T^{-1}x)e_p(-\xi x)=p^{-d}\sum_xf(x)e_p(-\xi Tx)=\hat{h}(T^*\xi),
	$$
	where $T^*$ is the conjugate operator to $T$, and thus $\|\hat{h}\|_1=\|\hat{f}\|_1$. Besides, since $\supp h =TA$, we see that (\ref{TA}) implies
	\begin{multline*}
	\|\hat{h}\|_1=p^{-d}\sum_{\xi_1,\ldots,\xi_d}\left|\sum_{x\in A'}h(a_x)e_p(-(y_2(x)\xi_2+\ldots+y_d(x)\xi_d))e_p(-x\xi_1)\right|\geq \\
	p^{-(d-1)}\sum_{\xi_2,\ldots,\xi_d}p^{-1}\sum_{\xi_1}  \left|\sum_{x\in A'}h(a_x)e_p(-(y_2(x)\xi_2+\ldots+y_d(x)\xi_d))e_p(-x\xi_1)\right|.
	\end{multline*} 
	For any fixed $\xi_2,\ldots,\xi_d$ we can apply the assumption of the theorem to the inner sum (because it is the Wiener norm of a function taking values in the set $E$). Since $|A'|=|A|$, each inner sum is at least $F(p,|A|p^{-1})$. 
	The claim follows.


\begin{thebibliography}{99}
			
\bibitem[TV]{TV} T.~Tao and V.~Vu, ``Additive Combinatorics'', Cambridge Stud. Adv. Math., Vol. 105.
		
		
\bibitem[Kon]{Kon} S.V.Konyagin, On a problem of Littlewood, Math. USSR-Izv., 18:2 (1982), 205–225. 
		
		
\bibitem[MPS]{MPS} McGehee, L.Pigno, B.Smith, Hardy's inequality and the $L_1$ norm of exponential sums, Annals of Math. 113 (1981), 613-618. 
		
\bibitem[GK]{GK} B.J.Green, S.V.Konyagin, On the Littlewood problem modulo prime, Canad. J.Math. 61 (2009), 141-164.
		
\bibitem[KS1]{KS small and big} 
		
S. V. Konyagin, I. D. Shkredov, A quantitative version of the Beurling-Helson theorem, Funct. Anal. Appl., 49:2 (2015), 110–121; arXiv: 1401.4429.
		
		
\bibitem[KS2]{KS medium} 
		
S. V. Konyagin, I. D. Shkredov, On the Wiener norm of subsets of Zp of medium size, J. Math. Sci., 218:5 (2016), 599–608.		

		
\bibitem[Sch]{Sch} T.Schoen, On the Littlewood conjecture modulo prime, Moscow Journal of Number Theory and Combinatorics, 1-5, 2016.
		
\bibitem[Teml]{Teml} V.N.Temlyakov, Multivariate Approximation, Cambridge University Press, 2018.
				
\bibitem[Sand1]{Sand1} T.Sanders, The Littlewood-Gowers problem, J.Anal.Math. 101 (2007), 123-162.
		
\bibitem[Sand2]{Sand2} T.Sanders, Bounds in Cohen's idempotent theorem, preprint, https://arxiv.org/pdf/1610.07092.pdf
		
		
	\end{thebibliography}
\end{document}